\documentstyle{article}
\begin{document}

\voffset=1cm
\newtheorem{prop}{Proposition}[section]
\newtheorem{th}[prop]{Theorem}
\newtheorem{cor}[prop]{Corollary}
\newtheorem{lm}[prop]{Lemma}
\newtheorem{df}[prop]{Definition}
\newtheorem{ex}[prop]{Example}
\newtheorem{remark}[prop]{Remark}
\newcommand{\bsquare}{\hbox{\rule{6pt}{6pt}}}\newcommand{\proof}[1]{\noindent{\bf Proof}\hspace{0.3cm}{#1}\hfill\bsquare \vspace{0.5cm}\par}
\renewcommand{\thefootnote}{\fnsymbol{footnote}}
\newcommand{\Ker}{{\rm{Ker}}~}
\renewcommand{\Im}{{\rm{Im}}~}
\newcommand{\Coker}{{\rm{Coker}}~}
\newcommand{\Hom}{{\rm{Hom}}}
\newcommand{\Supp}{\rm{Supp}}
\renewcommand{\H}{{\rm{H}}}
\renewcommand{\L}{{\rm{L}}}
\newcommand{\Res}{{\rm{Res}}}
\newcommand{\rank}{{\rm{rank}}}
\newcommand{\Pic}{{\rm{Pic}}}
\newcommand{\Div}{{\rm{Div}}}
\newcommand{\Der}{{\rm{Der}}}
\renewcommand{\O}{{\cal{O}}}
\newcommand{\U}{{\cal{U}}}
\newcommand{\G}{{\bf{G}}}

\title{Lefschetz pencils on a certain hypersurface in positive characteristic}
\author{Toshiyuki Katsura\thanks{Partially supported by JSPS Grant-in-Aid (C), No 24540053} 
}
\date{}
\maketitle

\rightline{\it Dedicated to Yujiro Kawamata on the occasion of  his 60th birthday}

\begin{abstract}
We examine Lefschetz pencils of a certain hypersurface in ${\bf P}^{3}$ over
an algebraically closed field of characteristic $p > 2$,
and determine the group structure of sections of the fiber spaces derived
from the pencils. 
Using the structure of a Lefschetz pencil, we give 
a geometric proof of the unirationality of Fermat surfaces 
of degree $p^a + 1$ with a positive integer $a$ 
which was first poved by Shioda \cite{S1}. 
As byproducts, we also see that on the hypersurface 
there exists a $(q^{3} + q^{2} + q + 1)_{q + 1}$-symmetric configuration
(resp. a $((q^{3} + 1)(q^{2} + 1)_{q + 1}, (q^{3} + 1)(q + 1)_{q^{2} + 1}$)-configuration)
 made up of the rational points over ${\bf F}_{q}$ (resp. 
over ${\bf F}_{q^{2}}$) and the lines over ${\bf F}_{q}$ 
(resp. over ${\bf F}_{q^{2}}$) with $q = p^{a}$. 
\end{abstract}

\section{Introduction}
Let $k$ be an algebaically closed field of charactersitic $p > 2$ and
we set  $q = p^{a}$ with a positive integer $a$.
Let $S$ be a hypersurface  in the 3-dimensional projective
space ${\bf P}^{3}$ defined by the equation
$x_{0}x_{1}^{q}-x_{1}x_{0}^{q}+x_{2}x_{3}^{q}-x_{3}x_{2}^{q}=0$.
The aim of this paper is to examine the structure of Lefschetz pencils 
on the surface $S$ and to  determine the singular fibers 
and sections of the fiber spaces derived from the pencils. 
In particular, in case of $p = q = 3$, this surface $S$ is a K3 surface. 
In fact, it is known that in this case the surface is a supersingular K3 surface
with Artin invariant 1. Our fiber space is a quasi-elliptic surface
with 10 singular fibers of type IV (for the existence of such a  quasi-elliptic surface,
 see H. Ito \cite{I}). 
As a corollary to our theory, we give a geometric proof of the fact 
that the Fermat surface 
of degree $q + 1$ is unirational, which was long ago proved by Shioda \cite{S1}
 (also see Rudakov-Shafarevich \cite{RS}). 

To examine Lefschetz pencils, we need to calculate rational
points and lines on $S$ defined over the finite field ${\bf F}_{q^2}$.
This part is known from various points of view (cf \cite{Seg}, \cite{SSL} and \cite{PT}),
but since we need to know the detailed structure to examine the Lefschetz pencils,
we give here a down-to-earth calculation for them.
Summing up our results, we conclude that on this surface there exists 
a $(q^{3} + q^{2} + q + 1)_{q + 1}$-symmetric configuration
(resp. a $((q^{3} + 1)(q^{2} + 1)_{q + 1}, (q^{3} + 1)(q + 1)_{q^{2} + 1}$)-configuration) 
made up of the rational points over ${\bf F}_{q}$ (resp. 
over ${\bf F}_{q^{2}}$) and the lines over ${\bf F}_{q}$ 
(resp. over ${\bf F}_{q^{2}}$) (also see \cite{PT} on the relation with 
the notion of finite generalized quadrangles).
In particular, in case of $p = q = 3$, we have a $(280_{4}, 112_{10})$-configuration 
on this K3 surface. 
Such a structure is related to the theory of Leech lattice and
these 112 lines correspond with Leech roots.
We examined the lattice structure of these lines in \cite{KK}.

The author would like to thank Professor Gerard van der Geer 
for suggesting him to examine the surface $S$ and for his advice.
The author would like to
thank Professors S. Kondo and T. Shioda for their valuable comments.
The author would also like to thank the referee for his careful reading and
useful comments. 
\section{Preliminaries}
We first recall the notion of a geometric realization of an abstract configuration.
A triple $\{{\cal A}, {\cal B}, R\}$, where  ${\cal A}, {\cal B}$ are non-empty
finite sets and $R \subset {\cal A}\times {\cal B}$ is a relation, is called an abstract
configuration if the cardinality of the set $R(x) = \{B \in {\cal B}~\mid~(x, B)\in R\}$
(resp. $R(B) = \{x \in {\cal A}~\mid~(x, B)\in R\}$) does not depend on $x \in {\cal A}$
(resp. $B \in {\cal B}$). Elements of ${\cal A}$
are called points, and elements of ${\cal B}$ are called blocks. Denoting 
by $\mid X\mid$ the number of elements in  a finite set $X$, we set
$$
v = \mid {\cal A}\mid,~b = \mid {\cal B}\mid, k =  \mid  R(x)\mid, r = \mid R(B)\mid.
$$
Then, the configuration is called a $(v_{k}, b_{r})$-configuration.
We have the relation $kv = br$.
Therefore, if $v = b$, then we have $k = r$.
In this case, the configuration is called a symmetric configuration.
Such a symmetric configuration is called $v_{k}$-configuration (for details, see
Dolgachev\cite{D}).

The most typical example of a geometric realization of an abstract configuration
is given by the projective plane over a finite field.
Let $p$ (resp. $a$)  be a prime number (resp. a positive integer)
and let ${\bf F}_{q}$ 
be a finite field with $q = p^{a}$ elements.
Then, in the projective plane ${\bf P}^{2}$ there are $q^2 + q + 1$ rational
points over  ${\bf F}_{q}$ and there are $q^2 + q + 1$ lines defined over
${\bf F}_{q}$. We see that $q + 1$ lines pass through each point, and on each line
there exist $q + 1$ points. We denote the set of these points by ${\cal A}$
and the set of these lines by ${\cal B}$. The relation $R$ consists of the pairs
of a point and a line which pass through the point. The triple ${\cal A}, {\cal B}, R$
gives a $(q^2 + q + 1)_{q + 1}$-symmetric configuration.

One more typical  configuration is given by Kummer surfaces.
Let $C$ be a non-singular complete curve of genus 2 
defined over an algebraically closed field
of characteristic $p \neq 2$. We consider the Jacobian variety $J(C)$.
Then, $C$ gives a principal polarization on $J(C)$, and by a suitable translation
we may assume that $C$ is invariant under the inversion $\iota$ of $J(C)$.
For a two-torsion point $a \in J(C)_{2}$, we denote by $T_{a}$ the translation
given by $a$. Then we have 16 curves $\{T_{a}C~\mid~a \in J(C)_{2}\}$.
We consider the quotient surface $J(C)/\langle \iota \rangle$, and
let $\pi  : J(C) \longrightarrow J(C)/\langle \iota \rangle$ be the projection.
Then, we have  the set ${\cal A}$ of 16 rational double points of type $A_{1}$ on $J(C)/\langle \iota \rangle$,
and we have the set ${\cal B} =\{ \pi(T_{a}C)\mid~a \in J(C)_{2}\}$
of 16 rational curves which are conics.
The relation $R$ consists of the pairs
of a point and a conic which pass through the point. The triple $\{{\cal A}, {\cal B}, R\}$
gives a $16_{6}$-symmetric configuration.

\section{Rational points over a finite field}
We consider the hypersurface $S$ in the 3-dimensional projective
space ${\bf P}^{3}$ which is defined by 
$$
(1) \quad \quad x_{0}x_{1}^{q}-x_{1}x_{0}^{q}+x_{2}x_{3}^{q}-x_{3}x_{2}^{q}=0
$$
It is easy to show that over ${\bf F}_{q^2}$ 
this surface is isomorphic to the Fermat surface defined by
$$
x_{0}^{q +1} +x_{1}^{q +1} +x_{2}^{q +1} +x_{3}^{q +1}  = 0.
$$
However, since the number of rational points over ${\bf F}_{q}$
of $S$ is different from the one of the Fermat surface, we see that $S$ is 
not isomorphic to the Fermat surface over ${\bf F}_{q}$. 
By the result in Weil \cite{W}, the number of ${\bf F}_{q^2}$-rational points
of the Fermat surface is known. Therefore, the number of 
${\bf F}_{q^2}$-rational points of $S$ is also known. 
However, to know the structure of the surface $S$ in detail
we give here a direct calculation 
of the number of ${\bf{F}}_{q^{2}}$-rational points.

Suppose $x_{0} \neq 0$.
To caluculate the rational points,  we may assume $x_{0}=1$. Then,  we have
the equation
$$
x_{1}^{q}-x_{1}=x_{3}x_{2}^{q}-x_{2}x_3^{q}.
$$
We have the following exact sequence of ${\bf{F}}_{q}$-vector spaces:
$$
(2)\quad \quad 0 \rightarrow {\bf{F}}_{q}\longrightarrow{\bf{F}}_{q^{2}} 
\stackrel{F-{\rm id}}{\longrightarrow} {\bf{F}}_{q^{2}}.
$$
Here, $F$ is the Frobenius morphism over ${\bf{F}}_{q}$
and ${\rm id}$ is the identity mapping.
We set 
$$
V=\{\alpha \in{\bf{F}}_{q^{2}}|\alpha^{q}=-\alpha\}.
$$
$V$ is a vector space over ${\bf{F}}_{q}$,
and we have
$$
{\rm Im} (F-{\rm id}) \subset V
$$
Since $\dim_{{\bf{F}}_{q}}V=\dim_{{\bf{F}}_{q}}{\rm Im}(F-{\rm id})=1$,  we see that
$V={\rm Im}(F-{\rm id})$.
     
Now, assume $x_{2}, x_{3} \in {\bf{F}}_{q^{2}}$. Then, $x_{2}^{q^2}=x_{2}$ and $x_{3}^{q2}=x_{3}$.
Therefore, we see $x_{3}x_{2}^{q}-x_{2}x_{3}^{q} \in  V$.
Hence, for each $x_{2}$ and $x_{3} \in {\bf{F}}_{q^{2}}$, we can find
$q$ numbers of $x_1$, using the exact sequence $(2)$ with $V={\rm Im}(F-{\rm id})$.
Hence, in this affine open set, the surface defined by the equation $(1)$ has 
$q \times q^{2} \times q^{2}=q^{5}$
rational points over ${\bf{F}}_{q^{2}}$.

Suppose now  $x_{0}=0$.
Then, the equation (1) becomes 
$$
x_{2}x_{3}^{q}-x_{3}x_{2}^{q}=0.
$$
Factorizing the left hand side, we have
$$
x_{2}x_{3}(x_{3}^{q-1}-x_{2}^{q-1})=x_{2}x_{3}\prod_{a \in {\bf{F}}_{q}^{*}} (x_{3}- ax_{2}).
$$
Here, ${\bf F}_{q}^{*}$ is the multiplicative group 
of non-zero elements of ${\bf F}_{q}$.
If $x_{2}=x_{3}=0$,  then we have  only one  rational point $(0, 1, 0, 0)$. 
 If $x_{2}=0$ and  $x_{3} \ne 0$, then the rational points are of the form 
$(0, *, 0, 1)$. Therefore, we have $q^{2}$ rational points.
If $x_{2} \ne 0$ and  $x_{3}=0$, then the rational points are of the form  
$(0, *, 1, 0) $. Therefore, we have    $q^{2}$ rational points.
If $x_{2} \ne 0$ and $x_{3} \ne 0$, then the rational points are of the form 
$(0, b , \alpha, a\alpha)$ with $b \in {\bf{F}}_{q^2}$, $a \in {\bf{F}}_{q}^{*}$, 
$\alpha \in {\bf{F}}_{q^2}^{*}$.
Moreover, if $b = 0$, the rational points are of the form $(0, 0, 1, a) $. Therefore, 
we have  $q-1$ rational points.
If $b \ne 0$, the rational points are of the form $(0, 1, \alpha, a\alpha)$.
Therefore, we have $(q-1)(q^{2}-1)$ rational points.

Hence, in total the number of rational points over ${\bf{F}}_{q^2}$ is equal to 
$$
\begin{array}{cc}
 q^{5}+1+q^{2}+q^{2}+(q-1)+(q-1)(q^{2}-1)  & =q^{5}+q^{3}+q^{2}+1\\
                 & =(q^{3}+1)(q^{2}+1) 
\end{array}
$$
Since the equation (1) contains all ${\bf{F}}_{q}$-rational points of ${\bf{P}}^{3}$, 
we see that the number of rational points over ${\bf{F}}_{q}$ is equal to 
$$ 
q^{3}+q^{2}+q+1      
$$

\section{Lines defined over a finite field} 
Now, we will count the number of lines defined over ${\bf F}_{q^2}$ 
(resp. ${\bf F}_{q}$) on the surface $S$. 
This number is already known in Tate\cite{T}, Segre\cite{Seg}, 
Sch\"utt-Shioda-van Luijk\cite{SSL} and Payne-Thas\cite{PT}, 
but to examine Lefschetz pencils we need to know 
how these lines sit on our surface $S$.

Suppose there exists a line $\ell$ defined over ${\bf{F}}_{q^2}$ 
(resp. ${\bf{F}}_{q}$)
on the surface (1). Then, on $\ell$ we have $q^{2}+1$ (resp. $q +1$)
rational points defined over ${\bf{F}}_{q^2}$ (resp. ${\bf{F}}_{q}$).
Therefore, any such line on (1) can be obtained by connecting two rational points
on $S$.

Take two rational points $P'=(\alpha_{0}, \alpha_{1}, \alpha_{2}, \alpha_{3})$, 
$Q' =(\beta_{0}, \beta_{1}, \beta_{2}, \beta_{3})$ on the surface (1) 
defined over ${\bf{F}}_{q^{2}}$(resp. ${\bf{F}}_{q}$), 
and assume that the line $\ell$ which connects $P'$ with $Q'$ 
lies on the surface (1). 
 Then,  for any $t \in k$,
 points $(\alpha_{0}+t\beta_{0}, \alpha_{1}+t\beta_{1}, \alpha_{2}+t\beta_{2}, \alpha_{3}+t\beta_{3})$   lie on
 the surface (1).  Substitute these points into (1).
 Since $P'$ and $Q'$ are points on the surface (1), we have
 $$
\alpha_{0} \beta_{1}^{q} t^{q}+\beta_{0}\alpha_{1}^{q}t - \beta_{0}^{q}\alpha_{1}t^{q} - \beta_{1}\alpha_{0}^{q}t
=\alpha_{3}\beta_{2}^{q}t^{q}+\beta_{3}\alpha_{2}^{q}t - \alpha_{2}\beta_{3}^{q}t^{q} -\beta_{2}\alpha_{3}^{q}t.
$$
Since  t is arbitrary, we have
$$
\begin{array}{l}
\alpha_{0}\beta_{1}^{q} - \alpha_{1}\beta_{0}^{q}=\alpha_{3}\beta_{2}^{q} - \alpha_{2}\beta_{3}^{q}, \\
\beta_{0}\alpha_{1}^{q} - \alpha_{0}^{q}\beta_{1}=\beta_{3}\alpha_{2}^{q} - \beta_{2}\alpha_{3}^{q}.
\end{array}
$$
These two equations have same solutions over ${\bf{F}}_{q^{2}}$
(resp ${\bf{F}}_{q}$).
Hence, the condition becomes
$$
(3) \quad \quad \quad  \alpha_{0}\beta_{1}^{q} - \alpha_{1}\beta_{0}^{q}=\alpha_{3}\beta_{2}^{q} - \alpha_{2}\beta_{3}^{q}
$$
Now, we consider the hyperplane $H'$ defined by 
$$
H':  \beta_{1}^{q}x_{0} - \beta_{0}^{q}x_{1}+ \beta_{3}^{q}x_{2} - \beta_{2}^{q}x_{3}=0
$$
This hyperplane is nothing but the tangent space of the surface (1) at the point $Q'$. 
By $(3)$, we see that $H'$ passes through the point $P'$.
Hence, any line defined over ${\bf{F}}_{q^{2}}$ (resp.${\bf{F}}_{q}$)
on the surface (1) is obtained as the lines cut by a tangent hyperplane at the rational points over ${\bf{F}}_{q^{2}}$ (resp.${\bf{F}}_{q}$).

Now, take a rational point $P=(\alpha, \beta, \gamma, \delta)$ on the surface (1) 
defined over ${\bf{F}}_{q^{2}}$ (resp. ${\bf{F}}_{q}$). 
Then, the tangent space $H$ of the surface (1) at $P$ is given by 
$$
(4)\quad \quad \quad  \beta^{q}x_{0} - \alpha^{q}x_{1}+ \delta^{q}x_{2}-\gamma^{q}x_{3}=0.
$$
Changing to inhomogeneous coordinates, without loss of generality 
we may assume the case $\gamma =1$.
Then, we have 
$$ 
x_{3}=\beta^{q}x_{0}-\alpha^{q}x_{1}+\delta^{q}x_{2}
$$
Substituting this into (1) and using 
$ -\alpha^{q}\beta+\beta^{q}\alpha=\delta^{q}-\delta$, 
we have an equation
$$
(x_{0}-\alpha x_{2})(x_{1}-\beta x_{2})
 \prod_{\epsilon\in {\bf{F}}_{q^{*}}}\{(x_{1}-\beta x_{2})- \epsilon(x_{0}-\alpha x_{2})\}=0.
$$
This means that the intersection of the surface (1) and the tangent space $H$
splits into $q+1$ lines defined over ${\bf{F}}_{q^2}$(resp. ${\bf{F}}_{q}$)
which intersect each other at the same point mutually transversely.
Since there exist $q^{2}+1$(resp. $ q+1$) rational points over ${\bf{F}}_{q^{2}}$
(resp. ${\bf{F}}_{q}$) on each line defined over 
${\bf{F}}_{q^{2}}$ (resp. ${\bf{F}}_{q}$), 
we conclude that on the surface (1) there are
$$
(q^{3}+1)(q^{2}+1)\times(q+1)\div(q^{2}+1)=(q^{3}+1)(q+1)
$$   
lines defined over ${\bf{F}}_{q^{2}}$ . We also see that on the surface (1) 
there exist
$$
(q^{3}+q^{2}+q+1)\times(q+1)\div(q+1)=q^{3}+q^{2}+q+1
$$
lines defined over ${\bf{F}}_{q}$. 

Hence, considering rational points and lines over ${\bf{F}}_{q^{2}}$ 
(resp.${\bf{F}}_{q}$) on the surface (1),
we have the following theorem.
\begin{th} On the hypersurface $S$ in ${\bf P}^{3}$ which is defined by 
$$
 x_{0}x_{1}^{q}-x_{1}x_{0}^{q}+x_{2}x_{3}^{q}-x_{3}x_{2}^{q}=0,
$$
there exist 
a $((q^{3} + 1)(q^{2} + 1)_{q + 1}, (q^{3} + 1)(q + 1)_{q^{2} + 1})$-configuration
and a $(q^{3} + q^{2} + q + 1)_{q + 1}$-symmetric configuration. 
\end{th}
\begin{remark}
In case $q = p = 3$, the surface $S$ given by $(1)$ is the supersingular
K3 surface with Artin invariant 1. In this case, our configuration is 
a $(280_{4}, 112_{10})$-configuration. We showed in $\cite{KK}$ that
$112$ lines correspond with Leech roots in the Picard lattice $Pic (S)$.
\end{remark}
\begin{remark}
In case $q = p$, the surface $S$ is related to the moduli space of supersingular 
K3 surfaces with Artin invariant $\sigma \leq 3$
(cf. Rudakov-Shafarevich $\cite{RS}$, p1520  and p1522, Theorem 2).
\end{remark}
\begin{remark}
In $\cite{S1}$, Shioda considered the hypersurface defined by
$x_{1}^{q}x_{2} + x_{1}x_{2}^{q} = x_{3}^{q}x_{0} + x_{3}x_{0}^{q}$
in ${\bf P}^{3}$. 
Over ${\bf F}_{q^{2}}$, this surface is isomorphic to the Fermat surface 
of degree $q + 1$ and also to our surface.
This surface is very similar to our surface $S$. However,
counting the number of rational points over  ${\bf F}_{q}$, 
we see that this surface is not isomorphic to our surface $S$ over ${\bf F}_{q}$.
\end{remark}
\begin{remark}
Let ${\cal A}$ and ${\cal B}$ be two sets, and $R$ be a relation between
${\cal A}$ and ${\cal B}$. The elements of ${\cal A}$ are called points
and the elements of ${\cal B}$ are called blocks. 
A triple $\{{\cal A}, {\cal B}, R\}$ is called a $t$-$(v, k, \lambda)$ design
if the following  three conditions hold.

$({\rm i})$  $\mid {\cal A}\mid = v$;

$({\rm ii})$  Every block $B \in {\cal B}$  relates to precisely $k$ points;

$({\rm iii})$ Every $t$ distinct points  together relates to precisely $\lambda$ blocks.

Using this notion, our 
$((q^{3} + 1)(q^{2} + 1)_{q + 1}, (q^{3} + 1)(q + 1)_{q^{2} + 1})$-configuration 
is a $1$-$((q^{3} + 1)(q^{2} + 1), q + 1, q+ 1)$ design.
\end{remark}
\begin{remark}
A $($finite$)$ generalized quadrangle  is an incidence structure $\{P, B, I \}$
in which P and B are disjoint nonempty sets, called points and lines, respectively,
and for which I is a symmetric point-line incidence relation which satisfies the following axioms:

$({\rm i})$ With an integer $t \geq 1$, each point is incident with $1 + t$ lines and two distinct points are incident with at most one line.

$({\rm ii})$ With an integer $s \geq 1$, each line is incident with $1 + s$ points and
two distinct lines are incident with at most one point.

$({\rm iii})$ If $x$ is a point and $L$ is a line not incident with $x$, then
there is a unique pair $(y, M) \in P \times B$ for which $ (x, M), (y, M), (y, L) \in I$.

The integers $s$ and $t$ are called the parameters of the generalized quadrangle
and $\{P, B, I\}$ is said to be order $(s, t)$ $($for the details, see $\cite{PT}$$)$.

Using this notion, 
our $((q^{3} + 1)(q^{2} + 1)_{q + 1}, (q^{3} + 1)(q + 1)_{q^{2} + 1})$-configuration 
is the generalized quadrangle of order $(q, q^2)$. From this point of view, 
this configuration is known in $\cite{PT}$, Chapter 3.
\end{remark}

\section{Lefschetz pencil}
On the surface $S$ defined by (1), we have $(q^{3} + 1)(q + 1)$ lines
defined over ${\bf F}_{q^2}$.
We take any line ${\ell}$ from these. 
Let $H$ and $H'$ be two different hyperplanes  such that $H \cap H' = \ell$.
Suppose that $H$ (resp. $H'$) is defined by the equation $L = 0$ (reps. $L' = 0$).
Then, our Lefschetz pencil on $S$ is defined as the pencil given 
by $\mu L + \mu'L' = 0$ with parameters $\mu, \mu'$. 
The line $\ell$ is the fixed component of the pencil.
Let $D + \ell$ be a general member of the pencil. As we explained in Section 4, 
by a suitable choice of $\mu$ and $\mu'$, we can find a member 
$\sum_{i = 1}^{q}\ell_{i} + \ell$ whose irreducible components $\ell_{i}$ $(i = 1, 2, \ldots, q)$ and $\ell$ are smooth lines which intersect
each other at the same point mutually transversely. Since $D + \ell$
is linearly equivalent to $\sum_{i = 1}^{q}\ell_{i} + \ell$, we have
$$
    (D + \ell, D) = (H, \sum_{i = 1}^{q}\ell_{i}) = q.
$$
On the other hand, we have
$$
   (D + \ell, D) = D^{2} + (\ell, D) = D^{2} + (\ell,  \sum_{i = 1}^{q}\ell_{i}) = D^2 + q.
$$
Therefore, we have $D^2 = 0$. Hence, our Lefschetz pencil gives rise to
a fiber space $f: S \longrightarrow {\bf P}^1$.
Here, one of general fibers coincides with $D$. 
We call this fiber space a Lefschetz fiber space.

Firstly, we consider the following special case.
\begin{lm}\label{lm:Lefschetz}
Let $\ell$ be a line on $S \subset {\bf P}^3$ given by $(1, 0, 0, s)$ 
with parameter $s$, and $f : S \longrightarrow {\bf P}^{1}$ be
 the Lefschetz fiber space by using the line $\ell$.
Then,  the general fiber is a rational curve with one singularity
and we have the singular fibers on the points $(t, 1) \in {\bf P}^1({\bf F}_{q^2})$.
\end{lm}
\proof{
Let $H$ (resp. $H'$) be the hyperplane defined by $x_{1} = 0$ (resp. $x_{2} = 0$).
Then, we have $H \cap H' = \ell$. The Lefschetz pencil is defined by
$$
   tx_{1} - x_{2} = 0,
$$
and the  Lefschetz fiber space is given by
$$
(5) \quad \quad \quad x_{0}x_{1}^{q-1}-x_{0}^{q}+tx_{3}^{q}-t^{q}x_{1}^{q-1}x_{3}=0
$$
with $t \in {\bf P}^{1}$. The cusp locus is given by $x_{1} = 0$
and the results follow from the equation (5).
}
\begin{th}\label{th} 
Let ${\bf F}_{q}$ be a finite field with $q = p^{a}$ elements.
Take any line $\ell$ on S and consider the Lefschetz fiber space
$f : S \longrightarrow {\bf P}^{1}$ with respect to $\ell$.
Then, the general fiber is a rational curve with one singular point and  
it has $q^2 + 1$ singular fibers. Each singular fiber consists of $q$ lines
which intersect each other at the same point mutually transversely.
\end{th}
\proof{
The general unitary group $GU_{4}(q)$ acts naturally on the surface $S$.
By the Witt theorem, we know that $GU_{4}(q)$ acts transitively
on the set of lines on S (cf. Appendix). 
This means that to show the first part of this theorem it  suffices to show it 
for a  line. Therefore, the first statement
follows from Lemma \ref{lm:Lefschetz}.

By the calculation of the previous section, the singular fibers exist over 
the ${\bf F}_{q^2}$-rational points of the base curve ${\bf P}^{1}$. 
Therefore, we have $q^2 + 1$ singular fibers. Again, by the calculation 
of the previous section, each singular fiber consists of $q$ nonsingular
rational curves which all intersect each other at the same point mutually
transversely. Therefore, we have in total $q\times (q^{2} + 1)$ lines 
in the singular fibers. The closure  of the singular loci of
general fibers is given by ${\ell}$. Therefore, it is a rational curve 
which is purely inseparable covering of degree
$q$ over the base curve.}
In the proof of Theorem \ref{th}, we call the closure ${\ell}$ of the singular loci 
of general fibers the cusp locus.
The following corollary was first proved by Shioda \cite{S1} 
(also see Rudakov-Shafarevich \cite{RS}). Our proof explains
the geometric meaning of the result.
\begin{cor}
The Fermat surface
$$
  x_{0}^{q + 1} + x_{1}^{q+1} + x_{2}^{q + 1} + x_{3}^{q + 1} = 0,
$$
 is unirational over an algebraically closed field in characteristic $p>0$.
\end{cor}
\proof{
The Fermat surface
$$
  x_{0}^{q + 1} + x_{1}^{q+1} + x_{2}^{q + 1} + x_{3}^{q + 1} = 0,
$$
is isomorphic to the surface $S$ over an algebraically closed field in characteristic $p>0$.
With the notation in Lemma \ref{lm:Lefschetz}, we consider the change of base
given by the Frobenius morphism $t = s^{q}$. Incidentally, this corresponds to 
the morphism from the singular locus to the base space which is given by the restriction
of the morphism $f$ to the singular locus. 
Then, by this change of base  we have a ruled surface 
over the projective line ${\bf P}^1$. 
Therefore, $S$ is unirational.
To show concretely by calculation,  first go to an inhomogeneous coordinate
with $x_{1} = 1$.
Then,  we have 
$$
x_{0} - x_{0}^{q} + s^{q}x_{3}^{q} - s^{q^{2}}x_{3} = 0.
$$
Setting $x_{0} - sx_{3} = y$, we have
$$
   (s - s^{q^2})x_{3} + y - y^{q} = 0,
$$
which shows $k(x_{0}, x_{3}, s) = k(s, y)$. Therefore, the surface $S$
is unirational.
}
The following lemma follows from a result on the representation of $GU_{4}(q)$
in Tate \cite{T}.
We give here a direct proof.
\begin{lm}\label{lm;line}
Any line on the surface S is defined over ${\bf F}_{q^2}$.
Any line on the Fermat hypersurface of degree $q + 1$ is also defined 
over ${\bf F}_{q^2}$.
\end{lm}
\proof{Take any line $\ell$ on S. 
Let $P = (\alpha_{0}, \alpha_{1}, \alpha_{2}, \alpha_{3})$
and $Q = (\beta_{0}, \beta_{1}, \beta_{2}, \beta_{3})$ be any two different 
points on $\ell$. To prove this lemma,
it suffices to find two different points on $\ell$ which are defined over ${\bf F}_{q^2}$.
With two parameters $s, t$, the point 
$$
(\alpha_{0}s+ \beta_{0}t, \alpha_{1}s + \beta_{1}t, \alpha_{2}s+\beta_{2}t, \alpha_{3}s + \beta_{3}t)
$$
exists on the surface $S$. Since $t$ and $s$ are arbitrary elements in $k$,  
we have 4 equations:
$$
\begin{array}{l}
  \alpha_{0}\alpha_{1}^q - \alpha_{1}\alpha_{0}^q + \alpha_{2}\alpha_{3}^q -\alpha_{3}\alpha_{2}^q = 0\\
  \beta_{0}\beta_{1}^q - \beta_{1}\beta_{0}^q + \beta_{2}\beta_{3}^q -\beta_{3}\beta_{2}^q = 0\\
  \alpha_{0}\beta_{1}^q - \alpha_{1}\beta_{0}^q + \alpha_{2}\beta_{3}^q -\alpha_{3}\beta_{2}^q = 0\\
  \beta_{0}\alpha_{1}^q - \beta_{1}\alpha_{0}^q + \beta_{2}\alpha_{3}^q -\beta_{3}\alpha_{2}^q = 0
\end{array}
$$
We consider 4-dimensional vector space $k^4$ and the following bilinear form on it:
$$
      u_{0}v_{0} + u_{1}v_{1} + u_{2}v_{2} + u_{3}v_{3}
$$
for $(u_{0}, u_{1}, u_{2}, u_{3}), (v_{0}, v_{1}, v_{2}, v_{3}) \in k^4$.
We consider the 2-dimensional subspace $V$ in $k^4$ generated by
$(\alpha_{1}^q, -\alpha_{0}^q, \alpha_{3}^q, -\alpha_{2}^q)$, 
$(\beta_{1}^q, -\beta_{0}^q, \beta_{3}^q, -\beta_{2}^q)$.
Then, considering the $q$-th powers of four equations above,
we see that 4 vectors 
$$
\begin{array}{l}
(\alpha_{0}, \alpha_{1}, \alpha_{2}, \alpha_{3})\\
(\beta_{0}, \beta_{1}, \beta_{2}, \beta_{3})\\
(\alpha_{0}^{q^2}, \alpha_{1}^{q^2}, \alpha_{3}^{q^2}, \alpha_{4}^{q^2})\\
(\beta_{0}^{q^2}, \beta_{1}^{q^2}, \beta_{2}^{q^2}, \beta_{3}^{q^2})
\end{array}
$$
are in the orthogonal subspace $V^{\perp}$ of $V$. Since $\dim V = 2$, we
have $\dim V^{\perp} = 2$.
Therefore, we have a relation
$$
{}^{t}\left(
\begin{array}{cccc}
\alpha_{0}^{q^{2}} & \alpha_{1}^{q^{2}} & \alpha_{2}^{q^{2}}&  \alpha_{3}^{q^{2}}\\
\beta_{0}^{q^{2}} & \beta_{1}^{q^{2}} & \beta_{2}^{q^{2}}&  \beta_{3}^{q^{2}}
\end{array}
\right)
= 
{}^{t}\left(
\begin{array}{cccc}
\alpha_{0} & \alpha_{1} & \alpha_{2}&  \alpha_{3}\\
\beta_{0} & \beta_{1} & \beta_{2}&  \beta_{3}
\end{array}
\right)
A
$$
with a $2 \times 2$-matrix $A$.  By the Lang-Steinberg theorem there exists a
regular $2 \times 2$-matrix $B$ such that $A = B^{-1}B^{(q^2)}$. 
Here, $B^{(q^2)}$ is the image of Frobenius map of degree $q^2$. 
Therefore, the first and the second rows of the matrix
$$
{}^{t}B^{-1}
\left(
\begin{array}{cccc}
\alpha_{0} & \alpha_{1} & \alpha_{2}&  \alpha_{3}\\
\beta_{0} & \beta_{1} & \beta_{2}&  \beta_{3}
\end{array}
\right)
$$
give two points on $\ell$  which are defined over ${\bf F}_{q^2}$.
}
Using the calculation over ${\bf F}_{q^2}$ in Section 4, we have the following result.
\begin{cor} 
The number of lines on $S$ is equal to $(q^3 + 1)(q + 1)$.
The number of lines on the Fermat hypersurface of degree $q + 1$ 
is also equal to $(q^3 + 1)(q + 1)$.
\end{cor}
\begin{th}
Under the same notation as in Theorem \ref{th},
the  group of the sections of the group scheme 
$S\setminus \ell \longrightarrow {\bf P}^{1}$ 
is isomorphic to $({\bf Z}/p{\bf Z})^{\oplus 4a}$.
\end{th}
\proof{
Take a section $C$ of 
the group scheme $S\setminus \ell \longrightarrow {\bf P}^{1}$.
Then, it intersects one of irreducible components 
of  each singular fiber with multiplicity one.
Therefore, it intersects a line $\ell'$ in each singular fiber with multiplicity one.
Since any singular fiber is given by an intersection of $S$ and a tangent space,
the section $C$ intersects the tangent space, which is 
a hyperplane in ${\bf P}^{3}$,
with multiplicity one. Hence, $C$ is a line on $S$. By Lemma \ref{lm;line},
the line on $S$ is defined over the finite field ${\bf F}_{q^2}$.
Therefore, by the consideration 
in the previous section, the hyperplane which is spanned by $C$ and $\ell'$
is a tangent space of $S$ at the intersection point of $C$ and $\ell'$, and
$C$ is one of $(q^{3} + 1)(q + 1)$ lines defined over  ${\bf F}_{q^2}$
which we already had. Since the number of singular fibrs is equal to $q^2 + 1$
and we have the cusp locus $\ell$ on $S$, we see that the number of
sections is equal to
$$
(q^{3} + 1)(q + 1) -q\times (q^{2} + 1) - 1 = q^4 = p^{4a}.
$$
Since the general fiber of $S\setminus \ell \longrightarrow {\bf P}^{1}$ is 
an additive group scheme ${\bf G}_{a}$ and any non-trivial torsion of 
${\bf G}_{a}$ is of order $p$, we know that these sections
form a group isomorphic to $({\bf Z}/p{\bf Z})^{\oplus 4a}$.
}

Finally,  we give a remark on a special case where the characteristic of the field k is 
equal to 3. Since it is known that the surface $S$ :
$$
x_{0}x_{1}^{3}-x_{1}x_{0}^{3}+x_{2}x_{3}^{3}-x_{3}x_{2}^{3}=0.
$$
is a supersingular K3 surface with Artin invariant 1, we summarize our results
in this interesting case.
By the consideration above,
we have 112 lines on $S$, which are all defined over ${\bf{F}}_{9}$.
Take any line among these 112 curves and make
the Lefschetz pencil $f : S \longrightarrow {\bf P}^{1}$ 
by using the line. Then, we have a quasi-elliptic fibration 
over the rational curve ${\bf P}^{1}$ with 10 singular fibers of type IV. 
We have just 10 ${\bf{F}}_{9}$-rational points on ${\bf P}^{1}$
on which the singular fibers lie. Hence, we have 30 lines in the singular fibers
and one line as the cusp locus which we use to make the Letschetz pencil.
The other lines are the sections of this quasi-elliptic surface. Therefore, we have the
following result.
\begin{th} Assume $q = p = 3$. 
Let $f : S \longrightarrow {\bf P}^{1}$ be the Lefschetz fiber space as above.
Then, it forms a quasi-elliptic surface with 10 singular fibers of type IV and
the Mordell-Weil group of this 
quasi-elliptic surface is isomorphic to $({\bf Z}/3{\bf Z})^{\oplus 4}$.
\end{th}
We note that the existence of quasi-elliptic surfaces with such singular fibers in characteristic 3 was shown by H. Ito. He also examined, in details,  
the structure of Mordell -Weil groups of quasi-elliptic surfaces (cf. Ito\cite{I}).

\section{Appendix}
We denote by $GL_{n}(q^{2})$ the general linear group which consists of 
all the regular $n\times n$-matrices with entries in ${\bf F}_{q^2}$.
For $x \in {\bf F}_{q^{2}}$, let  $x \mapsto \bar{x} = x^{q}$ be the automorphism
of ${\bf F}_{q^{2}}$ whose fixed field is $ {\bf F}_{q}$. We consider
the non-singular Hermitian form given by
$$
   x_{1}\bar{x_{3}} + x_{3}\bar{x_{1}} + x_{2}\bar{x_{4}} + x_{4}\bar{x_{2}}.
$$
The general unitary group $GU_{4}(q)$ is the subgroup of all elements of 
$GL_{4}(q^{2})$ that fix the non-singular Hermitian form.
We consider the hypersurface $S'$ defined by
$$
x_{1}\bar{x_{3}} + x_{3}\bar{x_{1}} + x_{2}\bar{x_{4}} + x_{4}\bar{x_{2}} = 0.
$$
in the 3-dimensional projective space ${\bf P}^{3}$. 
It is clear that $S'$ is isomorphic to the surface $S$ defined by the equation (1)
and $GU_{4}(q)$  acts on $S'$.
The following proposition is known, but for readers' convenience
we give here a concrete calculation. 
Since the order of $GU_{4}(q)$ is equal to
$$
(q + 1)q^{6}(q^{4} - 1)(q^3 + 1)(q^2 -1).
$$
(cf. \cite{CCNPW})  and the number of lines on $S$ is 
equal to $(q^{3} + 1)(q + 1)$,
we have, by the following proposition, an elementary proof
of the Witt theorem which we used in Section 5.
\begin{prop}\label{prop:fix}
Let $\ell$ be the line defined by $x_{1} = x_{2}= 0$. 
The order of the stabilizer of $GU_{4}(q)$ at $\ell$ is equal to
$$
q^{6}(q^{4} - 1)(q^{2} - 1).
$$
\end{prop}
We denote by $M_{2}(q^2)$ the set of all the $2\times 2$-matrices 
with entries in ${\bf F}_{q^2}$, and we first show the following lemma.
\begin{lm}\label{lm,GL}
We set $M = \{X \in M_{2}(q^2) ~\mid ~ {}^t\bar{X} = - X\}$.
Then we have $\mid M \mid = q^{4}$.
\end{lm}
\proof{
We set 
$$
X =
\left(
\begin{array}{cc}
   a & b \\
    c & d
\end{array}
\right).
$$
SInce $ {}^t\bar{X} = - X$, we have
$$
   a = - a^q , b = -c^q , c = -b^q, d = -d^q.
$$
The number of solutions  of $a = - a^q$ (resp. $d = -d^q$) in ${\bf F}_{q^2}$
 is equal to $q$, and the number of common solutions
 of $b = - c^q$ and  $c = -b^q$  in ${\bf F}_{q^2}$ is equal to $q^2$.
Hence, we have $\mid M \mid = q^{4}$.
}
Now, we prove Proposition \ref{prop:fix}.
Since the general unitary group $GU_{4}(q)$ fixes the  Hermitian form 
$x_{1}\bar{x_{3}} + x_{3}\bar{x_{1}} + x_{2}\bar{x_{4}} + x_{4}\bar{x_{2}}$,
the element $A \in GU_{4}(q)$ satisfies
$$
         AJ{}^t\bar{A} = J,
$$
where
$$
J =
\left(
\begin{array}{cc}
          0     &   E \\
         E     &   0
\end{array}
\right)
$$
with $2 \times 2$-identity matrix $E$.
Setting 
$$
 A =
\left(
\begin{array}{cc}
          A_{1}     &   A_{2}\\
          A_{3}     &   A_{4}
\end{array}
\right),
$$
with  $2 \times 2$-matrix $A_{i}$ ($i = 1, 2, 3, 4$), we have
$$
\begin{array}{l}
   A_{1}{}^t\bar{A_{2}} + A_{2}{}^t\bar{A_{1}} = 0, \\
  A_{1}{}^t\bar{A_{4}} + A_{2}{}^t\bar{A_{3}} = E,\\
  A_{3}{}^t\bar{A_{4}} + A_{4}{}^t\bar{A_{3}} = 0.
\end{array}
$$
Assume $A$ fixes the line $\ell$. This means that $A_{2} = 0$.
Therefore, we have
$$
  A_{1}{}^t\bar{A_{4}} = E,
A_{3}{}^t\bar{A_{4}} + A_{4}{}^t\bar{A_{3}} = 0.
$$
Therefore, $A_{4} \in GL_{2}(q^2)$, and the number of such matrices
is equal to $(q^{4} -  1)(q^{4} - q^{2})$. Since
$$
{}^t\overline{(A_{3}{}^t\bar{A_{4}})} = - A_{3}{}^t\bar{A_{4}},
$$
for each $A_{4} \in GL_{2}(q^2)$ we have, by Lemma \ref{lm,GL}, $q^4$ matrices 
in $M_{2}(q^2)$ which satisfy this equation. Hence, we conclude
that the order of the stabilizer at the line $\ell$ is equal to 
$q^{6}(q^{4} - 1)(q^{2} - 1)$.

\vspace{0.5cm}
\noindent
T.\ Katsura: Faculty of Science and Engineering, Hosei University,
Koganei-shi, Tokyo 184-8584, Japan

\noindent
E-mail address: toshiyuki.katsura.tk@hosei.ac.jp

\end{document}